\documentclass[review]{elsarticle}

\usepackage{lineno,hyperref}
\modulolinenumbers[5]

\journal{Journal of \LaTeX\ Templates}

\newcommand{\be}{\begin{equation}}
\newcommand{\ee}{\end{equation}}
\newcommand{\beq}{\begin{eqnarray}}
\newcommand{\eeq}{\end{eqnarray}}
\newcommand{\nbeq}{\begin{eqnarray*}}
\newcommand{\neeq}{\end{eqnarray*}}

\begin{document}

\begin{frontmatter}

\title{A Characterization of Exponential Distribution and the Sukhatme-R\'{e}nyi Decomposition of Exponential Maxima}

\author{George P. Yanev  and Santanu Chakraborty}
\address{School of Mathematical and Statistical Sciences \\
    The University of Texas Rio Grande Valley}

\begin{abstract}
A new characterization of the exponential distribution is established. It is proven that the well-known Sukhatme-Renyi
necessary condition is also sufficient for exponentiality. A method of proof due to Arnold and Villasenor based on the Maclaurin
series expansion of the density is utilized.
\end{abstract}

\begin{keyword}
characterization\sep exponential distribution\sep Sukhatme -- Renyi decomposition \sep maxima \sep random shifts
\MSC[2010] 62G30\sep  62E10
\end{keyword}

\end{frontmatter}

\linenumbers

\vspace{-0.3cm}\section{Introduction and Main Results}
\label{intro}
A number of characterizations of the exponential distribution are based on the distributional equation $X+T\stackrel{d}{=}Y$
involving a pair of random
variables $(X,Y)$ and a random translator (shift) variable $T$, independent of $X$.
Characterizations making use of this equation when $X$, $Y$, and $T$ are either order statistics or record
values were obtained in  Wesolowski and Ahsanullah (2004), Castano-Martinez et al. (2012), and Shah et al. (2014)
among others. In all studies so far the translator $T$
was assumed to follow a certain distribution. This restriction is removed in our theorem below.

Suppose $X_1, X_2,\ldots, X_n$ is a random sample of size $n\ge 2$ from a parent $X$ with absolutely continuous cdf $F$, such that $F(0)=0$.
Denote the maximum order statistic by $X_{n:n}$.

Arnold and Villasenor (2013) obtained a series of characterizations of the exponential distribution
based on a random sample of size two.  In particular, they proved that, under some additional conditions on the cdf $F$,
\[
X_{1}+\frac{1}{2}X_{2}\stackrel{{d}}{=}X_{2:2}
\]
characterizes the exponential distribution with some positive parameter.
They also made conjectures for extensions to larger sample sizes. In  Chakraborty and Yanev (2013) and Yanev and Chakraborty (2013) some of the results from Arnold and Villasenor (2013) were generalized to random samples of size $n\ge 3$.
For instance, it was proven in  Chakraborty and Yanev (2013), under the same assumptions on the cdf $F$ as in the case $n=2$, that for a fixed $n\ge 2$
\be \label{cy13}
X_{n-1:n-1}+\frac{1}{n}X_{n}\stackrel{{d}}{=}X_{n:n}
\ee
characterizes the exponential distribution.

The contribution of the present paper is twofold. (i) The characterization equation (\ref{cy13}) is extended to the case of maxima of $n$ and $n-s$ random variables for $1\le s\le n$.
(ii) The technique of proof from Arnold and Villasenor (2013) for a random sample of size two is expanded to the case of sample size $n\ge 2$ for any fixed $n$. The proof of the main result makes use of combinatorial identities, which might be of independent interest.  We believe that this technique will be useful in obtaining other characterization results in the future.

{\bf Theorem.}
 Let $X$ be a non-negative random variable with pdf $f(x)$. Assume that $f(x)$ is complex analytic for every $x$ and $f(0)> 0$.
 Let $n$ and $s$ be fixed integers such that $1\le s\le n-1$. If
\be \label{eqn1}
X_{n-s:n-s}+\frac{1}{n-s+1}X_{n-s+1}+\ldots +\frac{1}{n}X_n\stackrel{d}{=} X_{n:n},
\ee
then $X$ is exponential with some positive parameter.

It is well-known (cf. Conway, (1978), p.35) that every complex analytic function is infinitely differentiable and, furthermore, has a power series expansion about each point of its domain.

Note that the Theorem has been
applied in constructing goodness-of-fit tests for exponential distribution in  Jovanovic et al. (2015) and Volkova (2015).

The following direct corollary of the Theorem is of its own interest.

{\bf Corollary.} Let $X$ be a non-negative random variable with pdf $f(x)$. Assume that $f(x)$ is complex analytic for
every $x$ and $f(0)> 0$. If for fixed $n$
\be \label{cor}
X_1+\frac{1}{2}X_2+\frac{1}{3}X_3+\ldots +\frac{1}{n}X_n
\stackrel{d}{=}X_{n:n},
\ee
then $X$ is exponential with some positive parameter.

Equation (\ref{cor}) is a particular case (for maxima) of the
well-known Sukhatme-R\'{e}nyi decomposition (cf. Arnold et al., 2008, p.73) of
the $k$th order statistic in a random sample $X_1, X_2,\ldots,X_n$ from
an exponential distribution.
It is known (cf. Arnold and Villasenor, 2013) that if (\ref{cor}) holds for every $n$,
then necessarily $X_1,X_2,\ldots, X_n$ have a common exponential distribution.
Under the assumptions of the Corollary,  for $X_1,X_2,\ldots, X_n$ to be exponential it is sufficient that
(\ref{cor}) holds for one fixed $n$ only.

In the next section we state three lemmas, to be used in the proof of the Theorem. The main steps in the proof of the Theorem are outlined in Section~3. Details of the proof of the Theorem are given in Section~4. Section~5 contains
the proofs of Lemmas 1 and 2. Concluding remarks are given in the last section.

\vspace{-0.3cm}\section{Preliminaries}
\label{sec:1}
Introduce, for all non-negative integers $n$ and $i$, and any real number $x$,
\be \label{defH}
H_{n,i}(x):=\sum_{j=0}^n (-1)^j {n \choose j} (x-j)^i.
\ee
We start with identities involving $H_{n,i}(x)$, which may be of independent interest.

{\bf Lemma 1} Let $s$ and $r$ be positive integers. Then
\be \label{rec_formula}
\hspace{-1.8cm}(i) \qquad \sum_{j=0}^{r-1}{r \choose j}H_{s-1,j}(s)=H_{s,r}(s+1).
\ee
 \be \label{L2_1}
\hspace{-0.7cm}(ii) \qquad \sum_{j=0}^{r-1} {r \choose j+1}H_{s,j}(s+1)=\frac{1}{s+1}H_{s+1,r}(s+2).
\ee
\be \label{L2_2}
(iii) \qquad \sum_{j=0}^{r-1} (s+2)^{r-1-j}H_{s,j}(s+1)=\frac{1}{s+1}H_{s+1,r}(s+2).
\ee

Define $G_m(x):=F^m(x)f(x)$ for $m\ge 1$; $G_0(x):=F(x)$. Assuming (\ref{der}), we calculate the derivatives of $G_m(x)$ at 0 for $m\ge 1$.

{\bf Lemma 2}\  Let $m\ge 1$ and $d$ be integers, such that $d\ge -m$.
Assume $F(0)=0$. In case $d$ is positive, also assume,
\be \label{der}
f^{(k)}(0)= \left[\frac{f'(0)}{f(0)}\right]^{k-1}f'(0) \qquad k=1,\ldots, d,
\ee
then
\be \label{main}
G_m^{(m+d)}(0)=
\left\{
  \begin{array}{ll}
  \left[\frac{f'(0)}{f(0)}\right]^{d}f^{m+1}(0)H_{m,m+d}(m+1) & \mbox{if} \quad d\ge 0;
 \\
    0 & \mbox{if} \quad -m\le d<0.
  \end{array}
\right.
\ee

\noindent The third lemma, extracted from the proof of Theorem~1 in Arnold and Villasenor (2013), plays a central role in the proof of the Theorem.

{\bf Lemma 3}\ Let $X$ be a non-negative random variable with pdf $f(x)$. Assume that $f(x)$ is complex analytic for every $x$ and $f(0)>0$. If
\be \label{lemma}
f^{(k)}(0)=\left[\frac{f'(0)}{f(0)}\right]^{k-1}f'(0), \qquad k=1,2,\ldots,
\ee
then $X$ is exponential with some positive parameter.

Note that the assumptions for analyticity of $f(x)$ and $f(0)>0$ are implicitly used in the proof of Lemma~3 given in Arnold and Villasenor (2013).

\section{Outline of the Main Steps in the Proof of the Theorem}
\label{outline}

The proof of the Theorem can be divided into four steps as follows.

{\bf Step 1:}\  Define $d_j:=n-j+1$ and
$y_j:=z-x_s-\sum_{k=1}^{j-1} x_k$ for $1\le j\le s$. Show that the equality in distribution (\ref{eqn1}) is equivalent to
\beq \label{step0}
\lefteqn{\hspace{-0.5cm}\int_{0}^{z}G_{n-s-1}(x_s)\int_0^{y_1}\ldots \int_0^{y_{s-1}}\left(\prod_{j=1}^{s-1}f\left(d_jx_{j}\right)\right)f\left(d_sy_{s}\right)
\, dx_{s-1}\ldots dx_1dx_s } \nonumber \\
    & & =
f(z)\int_0^z \int_0^{x_{1}}\ldots \int_0^{x_{s-1}}\left(\prod_{k=1}^{s-1}f(x_k)\right)G_{n-s-1}(x_{s})\, dx_s\ldots dx_1.
\eeq

{\bf Step 2:}\
Denote
\[
r_j(t):=n-s+t+1-\sum_{k=1}^{j-1}i_k \qquad 1\le j\le s,\qquad t\ge 1,
\]
where $i_k$ are integers. We shall write $r_j$ instead of $r_j(t)$. Also, introduce
 \[
 a_{i_1,\ldots,i_s}  :=
   {d_s}^{r_s-i_s-1}\prod_{j=1}^{s-1}{d_j}^{i_j},\quad \quad
b_{i_1,\ldots,i_s}:={r_s\choose i_s+1}\prod_{j=1}^{s-1}{r_j+s-j\choose i_j}.
 \]
Prove (by differentiating (\ref{step0}) $(n+t)$ times with respect to $z$ and setting $z=0$) that (\ref{step0}) implies
\beq\label{newstep2}
\lefteqn{\sum_{i_1=0}^{r_1}\cdots \sum_{i_{s-1}=0}^{r_{s-1}}
\sum_{i_s=0}^{r_s-1}a_{i_1,\ldots,i_s}\left(\prod_{j=1}^{s-1}f^{(i_j)}(0)\right)
f^{(r_s-i_s-1)}(0)G^{(i_s)}_{n-s-1}(0)} \\
 & & \hspace{-1cm} =
\sum_{i_1=0}^{r_1}\cdots \sum_{i_{s-1}=0}^{r_{s-1}}
  \sum_{i_s=0}^{r_{s}-1} b_{i_1,\ldots,i_s}\left(\prod_{j=1}^{s-1}f^{(i_j)}(0)\right)f^{(r_s-i_s-1)}(0)G^{(i_s)}_{n-s-1}(0).\nonumber
\eeq

{\bf Step 3:}\
 Using Lemma~1, prove that
\beq \label{feqn}
\lefteqn{ \hspace{-2cm}\sum_{i_1=0}^{r_1}\cdots \sum_{i_{s-1}=0}^{r_{s-1}}\sum_{i_s=0}^{r_s-1}
a_{i_1,\ldots,i_s} H_{n-s-1, i_s}(n-s)} \\
& &  \hspace{-0.5cm} =\sum_{i_1=0}^{r_1}\cdots \sum_{i_{s-1}=0}^{r_{s-1}}\sum_{i_s=0}^{r_{s}-1}
b_{i_1,\ldots,i_s} H_{n-s-1,i_s}(n-s). \nonumber
 \eeq

{\bf Step 4:}\
Prove (\ref{lemma}) by induction using the results from Step~2 and Step~3.

The statement of the Theorem follows by Step~4 and Lemma~3.

\section{Proofs of the Steps in Section~3}
\label{sec:2}

Let $F_n(x)$ and $f_n(x)$ denote the cdf and pdf, respectively, of the maximum $X_{n:n}$.
Obviously, $F_{n}(x)=F^n(x)$.

\subsection{Proof of Step 1}
\label{subsec:5.1}
Let $f_{n-1,n}(x)$ denote the density of
$X_{n-1}/(n-1)+X_{n}/n$.
Setting $s=2$ in (\ref{eqn1}), for the density $f_{LHS}(z\ |s=2)$, say, of the left-hand side of (\ref{eqn1}) we find
\nbeq
\lefteqn{f_{LHS}(z \ | s=2) = \int_{0}^{z}f_{n-2}(x_2)f_{n-1,n}(z-x_2)\, dx_2} \\
& = & \int_{0}^{z}\!(n-2)G_{n-3}(x_2)n(n-1)\!\int_{0}^{z-x_2}\!f(nx_1)f((n-1)(z-x_2-x_1))\, dx_1dx_2 \\
& = & (n)_{3}\int_{0}^{z}G_{n-3}(x_2)\int_{0}^{z-x_2}\!f(nx_1)f((n-1)(z-x_2-x_{1}))\, dx_1dx_2,
\neeq
where $(n)_3:=n(n-1)(n-2)$. Setting $s=3$ in (\ref{eqn1}),
we have
\nbeq
\lefteqn{f_{LHS}(z \ | s=3)=  (n)_4\int_{0}^{z}\!G_{n-4}(x_3) \!\int_{0}^{z-x_3}\!\int_0^{z-x_3-x_1}\!\int_0^{z-x_3-x_1-x_2}}\\
& & \times f(nx_1)f((n-1)x_2)f((n-2)(z-x_3-x_1-x_2))\, dx_2dx_1dx_3.
\neeq
Repeating this argument, we obtain
 for any $s$ such that $2\le s\le n-1$,
\beq \label{LHS1}
\lefteqn{\frac{f_{LHS}(z)}{(n)_{s+1}}}\\
& & \hspace{-1cm}=\int_{0}^{z}\! \!G_{n-s-1}(x_s)\int_0^{y_1}\! \!\ldots \int_0^{y_{s-1}}\! \!\left(\prod_{j=1}^{s-1}f\left(d_jx_{j}\right)\right)f\left(d_sy_{s}\right)\, dx_{s-1}\ldots dx_1dx_s. \nonumber
\eeq
For the density $f_{RHS}(z)$, say, of the right-hand side of (\ref{eqn1}),
we have
\nbeq
\lefteqn{\hspace{-2cm} f_{RHS}(z)  =  n f(z)F_{n-1}(z)} \\
 & & \hspace{-2cm} =  nf(z)\int_0^z f_{n-1}(x_1)\, dx_1   \nonumber \\
 & & \hspace{-2cm} = n(n-1)f(z)\int_0^z f(x_1)F_{n-2}(x_1)\, dx_1   \nonumber \\
    & & \hspace{-2cm} =n(n-1)(n-2)f(z)\int_0^z\!\int_0^{x_1}f(x_1)f(x_2) F_{n-3}(x_2)\, dx_2dx_1. \nonumber
 \neeq
Repeating this argument $(s-2)$ more times we obtain
\be  \label{RHS1}
\hspace{-0.2cm}\frac{f_{RHS}(z)}{(n)_{s+1}}
=f(z)\int_0^z \int_0^{x_{1}}\ldots \int_0^{x_{s-1}}\left(\prod_{k=1}^{s-1}f(x_k)\right)G_{n-s-1}(x_{s})\, dx_s\ldots dx_1.
 \ee
Combining (\ref{LHS1}) and  (\ref{RHS1}) we obtain (\ref{step0}).

\subsection{Proof of Step 2}
Define
\[
K_{n,s-1}(y_1):=\int_0^{y_1}\! \!\ldots \int_0^{y_{s-1}}\! \!\left(\prod_{j=1}^{s-1}f\left(d_jx_{j}\right)\right)f\left(d_sy_{s}\right)\, dx_{s-1}\ldots dx_1.
\]
Observing that $K_{n,s-1}^{(i)}(0)=0$ when $i<s-1$ and  $G_{d_{s+2}}^{(i)}(0)=0$ for $i<d_{s+2}$, for the $(n+t)$th derivative of the left-hand side of (\ref{step0}) at 0, we obtain
\beq \label{n_t_der}
\lefteqn{\frac{d}{dz^{n+t}}\left\{\int_{0}^{z}G_{d_{s+2}}(x_s)K_{n,s-1}(z-x_s) dx_s\right\}\mid_{z=0}}\\
& & \hspace{-0.7cm} =
\frac{d}{dz^{n+t-1}}\left\{G_{d_{s+2}}(z)K_{n,\, s-1}(0)+\int_{0}^{z}G_{d_{s+2}}(x_s)K'_{n,s-1}(z-x_s) dx_s\right\}\mid_{z=0} \nonumber\\
& & \hspace{-0.7cm} =
\frac{d}{dz^{n+t-s}}\left\{G_{d_{s+2}}(z)K_{n,s-1}^{(s-1)}(0)+\int_{0}^{z}G_{d_{s+2}}(x_s)K_{n,s-1}^{(s)}(z-x_s)\, dx_s\right\}\mid_{z=0}\nonumber\\
& & \hspace{-0.7cm} =
\sum_{i=n-s-1}^{n-s+t}G_{d_{s+2}}^{(i)}(0)K_{n,s-1}^{(n+t-1-i)}(0).\nonumber
\eeq
Using the recursive relation
\[
K_{n,s-1}(u) =  \int_{0}^{u}f(nx)K_{n-1,s-2}(u-x)\, dx\qquad 3\le s\le n-1,
\]
one can show by induction that the $m$th derivative of $K_{n,s-1}(u)$ at 0 for
any $m\ge s-1$ and any fixed $n\ge 2$ is given by
\begin{equation} \label{eq:16}
K_{n,s-1}^{(m)}(0)=\sum_{i_{1}=0}^{l_1}\cdots\sum_{i_{s-1}=0}^{l_{s-1}}\left(\prod_{j=1}^{s-1}d_{j}^{i_{j}} f^{(i_{j})}(0)\right)d_{s}^{l_{s}}f^{(l_{s})}(0),
\end{equation}
where $l_j=m-s+1-\sum_{l=1}^{j-1}i_l$ for $1\le j\le s$. We omit the derivation of (\ref{eq:16}) here.
Substituting (\ref{eq:16}) into (\ref{n_t_der}) and changing the indexes of summation, one can see that the last sum in (\ref{n_t_der}) equals
\be \label{new_ls}
\sum_{i_1=0}^{r_1}\cdots \sum_{i_{s-1}=0}^{r_{s-1}}\left(\prod_{j=1}^{s-1}d_j^{i_j}f^{(i_j)}(0)\right)
\sum_{i_s=0}^{r_s-1}
d_s^{r_s-i_s-1}f^{(r_s-i_s-1)}(0)G^{(i_s)}_{d_{s+2}}(0).
\ee
where, as before, $r_{j}=n-s+t+1-\sum_{k=0}^{j-1}i_{k}$ for $1\le j\le s$.
Thus, we have obtained the left-hand side of (\ref{newstep2}).

We turn to the right-hand side of (\ref{newstep2}). Denote
\[
L(x_1|x_2,\ldots,x_s):=\int_0^{x_1}\ldots \int_0^{x_{s-1}}\left(\prod_{k=1}^{s-1}f(x_k)\right)G_{d_{s+2}}(x_{s})\, dx_s\ldots dx_2.
\]
With this notation for the $(n+t)$th derivative of $f_{RHS}(z)/(n)_{s+1}$ at 0, we find
\nbeq \label{RHS2}
\lefteqn{\frac{d}{dz^{n+t}}
\left\{f(z)\int_0^z L(x_1|x_2,\ldots,x_s)\, dx_1
\right\}\mid_{z=0}} \\
 & & = \sum_{i_1=0}^{n+t} {n+t \choose i_1} f^{(i_1)}(0)
\frac{d}{dz^{n+t-i_1}}\left\{
\int_0^z f(x_1)L(x_1|x_2,\ldots,x_s)\, dx_1
\right\}\mid_{z=0}  \nonumber \\
 & & = \sum_{i_1=0}^{n+t-1} {n+t \choose i_1} f^{(i_1)}(0)\frac{d}{dz^{n+t-1-i_1}}\left\{
f(z)\int_0^z L(x_2|x_3,\ldots,x_s)\, dx_2
\right\}\mid_{z=0} \nonumber
 \neeq
Recall that $r_j:=n+t-s+1-\sum_{k=1}^{j-1} i_k$ for $j=1,\ldots,s$. Repeating the last argument, it is not difficult to obtain
 \beq \label{new4}
 \lefteqn{\frac{d}{dz^{n+t}}
\left\{f(z)\int_0^z L(x_1|x_2,\ldots,x_s)\, dx_1
\right\}\mid_{z=0}}\\
& & \hspace{-1cm} = \sum_{i_1=0}^{r_1}\cdots \sum_{i_{s-1}=0}^{r_{s-1}}
\left(\prod_{j=1}^{s-1}{r_j+s-j \choose i_j}f^{(i_j)}(0)\right)\frac{d}{dz^{r_s}}
\left\{f(z)\int_0^z G_{d_{s+2}}(x_s)dx_s\right\}\mid_{z=0}\nonumber  \\
& & \hspace{-1cm}=
 \sum_{i_1=0}^{r_1}\cdots \sum_{i_{s-1}=0}^{r_{s-1}}
\left(\prod_{j=1}^{s-1}{r_j+s-j \choose i_j}f^{(i_j)}(0)\right)
\sum_{i_s=0}^{r_{s}-1}{r_{s}\choose i_{s}+1}
f^{(r_s-i_s-1)} (0)G^{(i_s)}_{d_{s+2}}(0). \nonumber
 \eeq
Combining (\ref{new_ls}) and  (\ref{new4}) we prove Step~2.

\subsection{Proof of Step 3}
 \label{subsec:5.4}

We shall simplify the right-hand side of (\ref{feqn}), working on the most inner sum first and moving to the outer ones later. Applying (\ref{L2_1}), we see that
\nbeq
\lefteqn{\hspace{-1cm}\sum_{i_{s-1}=0}^{r_{s-1}}{r_s+1\choose i_{s-1}}\sum_{i_s=0}^{r_s-1}{r_s\choose i_s+1}H_{n-s-1,i_s}(n-s)}\\
   & = &
\frac{1}{n-s}\sum_{i_{s-1}=0}^{r_{s-1}}{r_s+1\choose i_{s-1}}H_{n-s,r_{s-1}-i_s}(n-s+1)\nonumber \\
   & = &
\frac{1}{(n-s)(n-s+1)}H_{n-s+1,r_{s-2}+1-i_{s-2}}(n-s+2).
\neeq
Furthermore, since $H_{n-s+1,0}(n-s+1)=0$, applying (\ref{L2_1}) again, we have
\nbeq
\lefteqn{\hspace{-1cm}\sum_{i_{s-2}=0}^{r_{s-2}}{r_{s-2}+2\choose i_{s-2}}H_{n-s+1,r_{s-2}+1-i_{s-2}}(n-s+2)}\\
  & = &
\sum_{i_{s-2}=0}^{r_{s-2}+1}{r_{s-2}+2\choose i_{s-2}}H_{n-s+1,r_{s-2}+1-i_{s-2}}(n-s+2)\\
   & = &
\frac{1}{(n-s+2)}H_{n-s+2,r_{s-3}+2-i_{s-3}}(n-s+3).
\neeq
Repeating this argument for the rest of the sums on the right-hand side  of (\ref{feqn}), we find
\beq \label{simple_rhs}
\lefteqn{\hspace{-4cm}\sum_{i_1=0}^{r_1}\cdots \sum_{i_{s-1}=0}^{r_{s-1}}\left(\prod_{j=1}^{s-1}{r_j+s-j\choose i_j}\right)\sum_{i_s=0}^{r_{s}-1}
{r_s\choose i_s+1}H_{n-s-1,i_s}(n-s)} \\
& & =\frac{H_{n-1,n+t}(n)}{(n-1)_s}, \nonumber
\eeq
where $d_{s+2}=n-s-1$. Let us turn to the left-hand side of (\ref{feqn}). Recall that $H_{t,i_s}(\cdot)=0$ for $0\leq i_s\leq t-1$.
Similarly to the arguments in the simplification of the right-hand side above, applying (\ref{L2_2}) instead of (\ref{L2_1}), we obtain
\be \label{simple_lhs}
\hspace{-0.5cm}\sum_{i_1=0}^{r_1}\cdots \sum_{i_{s-1}=0}^{r_{s-1}}\left(\prod_{j=1}^{s-1}d_j^{i_j}\right)\sum_{i_s=0}^{r_s-1}
d_s^{r_s-i_s-1} H_{n-s-1, i_s}(n-s)
 =  \frac{H_{n-1,n+t}(n)}{(n-1)_s}.
\ee
Equations (\ref{simple_rhs}) and (\ref{simple_lhs}) imply (\ref{feqn}), which completes the proof of Step~3.

\subsection{Proof of Step 4}
 \label{subsec:5.3}
Denote $c_{i_1,\ldots,i_s}:=a_{i_1,\ldots,i_s}-b_{i_1,\ldots,i_s}$.
With this notation and taking into account that $G^{(i_s)}_{n-s-1}(0)=0$ when $i_s< n-s-1$, we write $(\ref{newstep2})$ as
\be \label{reduced_sum}
\hspace{-0.3cm} \sum_{i_1=0}^{r_1}\cdots \sum_{i_{s-1}=0}^{r_{s-1}} \sum_{i_s=n-s-1}^{r_{s}-1}c_{i_1,\ldots,i_s}\left(\prod_{j=1}^{s-1}f^{(i_j)}(0)\right)f^{(r_s-i_s-1)}(0)G^{(i_s)}_{n-s-1}(0)=0.
\ee
We shall prove (\ref{lemma}) by (strong) induction on $k$. The base case $k=1$ is trivial. Assuming (\ref{lemma}) for $k\le t$, we shall prove it for $k=t+1$, where
$t$ stands for any positive integer.
First, observe that since the order of the derivative of $f(x)$ in (\ref{newstep2}) must be nonnegative, we have $r_s-i_s-1\ge 0$. Combining this with $i_s\ge n-s-1$, we see that
\be \label{index_ineq}
\sum_{k=1}^{s-1}i_k \le t+1.
\ee
To extract the terms with a factor $f^{(t+1)}(0)$, we shall split the sum in (\ref{reduced_sum}) into two as follows
\beq  \label{RHS9}
\lefteqn{\hspace{-1.8cm}\sum_{{\cal I}\setminus {\cal I}_0}c_{i_1,\ldots,i_s}\left(\prod_{j=1}^{s-1}f^{(i_j)}(0)\right) f^{(r_s-i_s-1)}(0)G^{(i_s)}_{n-s-1}(0)   }\\
 & &  + f^{s-1}(0)f^{(t+1)}(0)G^{(n-s-1)}_{n-s-1}(0)\sum_{{\cal I}_0}c_{i_1,\ldots,i_s}=0, \nonumber
 \eeq
where ${\cal I}=\{(i_1,\ldots ,i_s):\ 0\le i_j\le r_j, \ 1\le j\le s-1, 1\le i_s\le r_s-1\}$ and ${{\cal I}_0}$ is the set of vectors $(i_1,\ldots,i_s)$ such that $i_s=n-s-1$ and among the first $s-1$ components: (i) all are zeros or (ii) exactly one is $t+1$ and the others are zeros. Notice that by Lemma~2
\be \label{G_n-s-1_term}
G^{(n-s-1)}_{n-s-1}(0)=f^{n-s}(0)H_{n-s-1,n-s-1}(n-s).
\ee
Consider the first sum in (\ref{RHS9}) (the one over ${\cal I}\setminus {\cal I}_0$). Inequality (\ref{index_ineq}) along with the definition of the index set ${\cal I}\setminus {\cal I}_0$ implies that all derivatives of $f(x)$ included in the product term have order
less than or equal to $t$. Therefore, applying the induction assumption to $f^{(i_j)}(0)$ for $i_j\ge 1$, we have
\be \label{product_term}
\prod_{j=1}^{s-1}f^{(i_j)}(0)=
\left\{
  \begin{array}{ll}
  f(0)[f'(0)]^{\sum_{k=1}^{s-1}i_k} & \mbox{if} \quad (i_1,\ldots,i_{s-1})\ne (0,\ldots,0);
 \\
    1 & \mbox{if} \quad (i_1,\ldots,i_{s-1})= (0,\ldots,0).
  \end{array}
\right.
\ee
It is not difficult to see that over the index set ${\cal I}\setminus {\cal I}_0$ we have $r_s-i_s-1\le t$ and therefore, applying the induction assumption, we obtain for $n-s-1\le i_s\le r_s-1$
\be \label{single_term}
 f^{(r_s-i_s-1)}(0)=\left[\frac{f'(0)}{f(0)}\right]^{r_s-i_s-2}f'(0).
\ee
It remains to study the factor $G^{(i_s)}_{n-s-1}(0)$. Since $i_s\le r_s-1\le n-s+t-\sum_{k=1}^{s-1}i_k$, we have that $i_s-(n-s-1)\le t+1-\sum_{k=1}^{s-1}i_k$.
We consider two cases as follows. (i) Let $\sum_{k=1}^{s-1}i_k\ge 1$. Then $i_s-(n-s-1)\le t$ and, under the induction assumption, applying Lemma~2 with $m=n-s-1$ and $d=i_s-(n-s-1)\le t$, we have
\be \label{G_term}
G_{n-s-1}^{(i_s)}(0)=\left[\frac{f'(0)}{f(0)}\right]^{i_s-n+s+1}f^{n-s}(0)H_{n-s-1,i_s}(n-s).
\ee
(ii) Let  $\sum_{k=1}^{s-1}i_k=0$. If $i_s\le n-s+t-1$, then (\ref{G_term}) holds. If $i_s=n-s+t$, then we see that
\be \label{c_term}
c_{0,\ldots,0,n-s+t}=0.
\ee
Combining (\ref{G_n-s-1_term})-(\ref{c_term}), it is not difficult to obtain that, under the induction assumption, (\ref{RHS9}) implies
\nbeq
\lefteqn{\hspace{-2cm}\left[\frac{f'(0)}{f(0)}\right]^{t}\! \! f'(0)
\sum_{{\cal I}\setminus {\cal I}_0} c_{i_1,\ldots,i_s} H_{n-s-1, i_s}(n-s)}\\
& &  + f^{(t+1)}(0)\sum_{{\cal I}_0}c_{i_1,\ldots,i_s}H_{n-s-1,n-s-1}(n-s)=0.
\neeq
Thus, to prove (\ref{lemma}) for $k=t+1$, it is sufficient to prove
\[
\sum_{{\cal I}\setminus {\cal I}_0} c_{i_1,\ldots,i_s} H_{n-s-1, i_s}(n-s)
 + \sum_{{\cal I}_0}c_{i_1,\ldots,i_s}H_{n-s-1,n-s-1}(n-s)=0
 \]
 or, equivalently,
 \[
 \sum_{\cal I}a_{i_1,\ldots,i_s} H_{n-s-1, i_s}(n-s)
 =\sum_{\cal I}b_{i_1,\ldots,i_s}H_{n-s-1,i_s}(n-s).
 \]
This is equivalent to (\ref{feqn}) proven to be true in Step~3. Therefore, the proof of Step~4 is complete.

\section{Proofs of Lemma 1 and Lemma 2}
\label{sec:3}
It is known (cf.  Ruiz, 1996) that for any non-negative integer $n$ and real $x$
\[
H_{n,i}(x)=
\left\{
  \begin{array}{ll}
    n! & \quad \mbox{if} \quad i=n; \\
    0 & \quad \mbox{if} \quad 0\le i< n.
  \end{array}
\right.
\]
This information will be useful in the proofs of the lemmas that follow.

{\bf Proof of Lemma 1.} (i) By the definition of $H_{n,i}(x)$ in (\ref{defH}), we obtain
\nbeq
\lefteqn{\sum_{j=0}^{r-1}{r \choose j}H_{s-1,j}(s)=\sum_{i=0}^{s-1}(-1)^i{s-1 \choose i}\sum_{j=0}^{r-1} {r \choose j}(s-i)^{j}}\\
    & & =\sum_{i=0}^{s-1}(-1)^i{s-1 \choose i}\left[(s+1-i)^r-(s-i)^r\right]\\
    & & =(s+1)^r-\left[s^r+{s-1 \choose 1}s^r\right]+ \ldots +(-1)^{s-1}\left[{s-1 \choose s-2}2^s+2^s\right]+(-1)^s\\
    & & =(s+1)^r - {s \choose 1}s^r +\ldots +(-1)^{s-1}{s \choose s-1}2^r+(-1)^s\\
    & & =\sum_{j=0}^s (-1)^j{s \choose j} (s+1-j)^r \\
    & & = H_{s,r}(s+1).
    \neeq
(ii) Indeed, using the definition of $H_{s,j}(x)$ in (\ref{defH}), we have
\nbeq
\lefteqn{\sum_{j=0}^{r-1}{r \choose j+1}H_{s,j}(s+1)=
 \sum_{j=0}^{r-1} {r \choose j+1} \sum_{i=0}^s (-1)^i {s \choose i} (s+1-i)^{j}}\\
 & = &
\sum_{i=0}^s (-1)^i {s \choose i} \sum_{k=1}^r {r \choose k}  (s+1-i)^{k-1}\\
& = &
\sum_{i=0}^s (-1)^i {s \choose i} \frac{1}{s+1-i}\left[\sum_{k=0}^r {r \choose k}  (s+1-i)^{k}-1\right]\\
& = &
\frac{1}{s+1}\sum_{i=0}^s (-1)^i {s+1 \choose i} [(s+2-i)^r-1]\\
& = &
\frac{1}{s+1}\sum_{i=0}^{s+1} (-1)^i {s+1 \choose i}(s+2-i)^r \\
& = &
\frac{1}{s+1}H_{s+1,r}(s+2).
\neeq
(iii) We have
\nbeq
\lefteqn{\hspace{-2cm}\sum_{j=0}^{r-1}(s+2)^{r-1-j}H_{s,j}(s+1)=
\sum_{j=0}^{r-1} (s+2)^{r-1-j} \sum_{i=0}^s (-1)^i {s \choose i} (s+1-i)^{j}}\\
 & = &
 \sum_{i=0}^s (-1)^i {s \choose i}(s+2)^{r-1} \sum_{j=0}^{r-1}
\left(\frac{s+1-i}{s+2}\right)^j\\
  & = &
 \sum_{i=0}^s (-1)^i {s \choose i}\frac{1}{i+1}\left[(s+2)^r-(s+1-i)^r\right]\\
  & = &
   \frac{1}{s+1}\sum_{i=0}^s (-1)^i{s+1 \choose i+1} \left[(s+2)^r-(s+1-i)^r\right]\\
  & = &
\frac{1}{s+1}\sum_{j=0}^{s+1} (-1)^j {s+1 \choose j} (s+2-j)^r \\
&  = &
\frac{1}{s+1}H_{s+1,r}(s+2).
\neeq

{\bf Proof of Lemma 2.}
(i) If $-m\le d<0$, then $G^{(m+d)}_m(0)=0$ because all the terms in the expansion of
$G^{(m+d)}_m(0)$ have a factor $F(0)=0$.

(ii) Let $d=0$. We shall prove (\ref{main}) by induction on $m$. One can verify directly the case $m=1$.
Assuming (\ref{main}) for $m=k$, we shall prove it for $m=k+1$.
Since $G_{k+1}(x)=F(x)G_{k}(x)$, applying (i), we see that
\nbeq
G^{(k+1)}_{k+1}(0) & = & \sum_{j=0}^{k+1} {k+1 \choose j}F^{(j)}(0) G^{(k+1-j)}_k(0)\\
    & & \hspace{-2cm}=  F(0)G^{(k+1)}_k(0)+(k+1)F'(0)G^{(k)}_k(0)+\sum_{j=2}^{k+1} {k+1 \choose j}F^{(j)}(0) G^{(k+1-j)}_k(0)\\
    & & \hspace{-2cm} = (k+1)!f^{k+2}(0),
\neeq
which completes the proof of (ii).

(iii) Let $d>0$ and $m$ be any positive integer.
 For simplicity, we will write $f^{(j)}:=f^{(j)}(0)$ below.

(a) Let $m=1$. If $d=1$, then we have $G_1^{(2)}(0)=3f'f=f'fH_{1,2}(2)$ since
$H_{1,2}(2)=3$.
Thus, (\ref{main}) is true for $d=1$. Next, assuming (\ref{main}) for $G_1^{(k)}(0)$, we shall prove it for $G_1^{(k+1)}(0)$.
Since $G_1(x)=F(x)f(x)$, using (\ref{der}) we obtain
\nbeq
G_1^{(k+1)}(0)
    & = & \sum_{j=1}^{k+1} {k+1 \choose j} f^{(j-1)}f^{(k+1-j)}\\
    & = & \sum_{j=1}^{k+1} {k+1 \choose j} \left(\frac{f'}{f}\right)^{j-2}f'\left(\frac{f'}{f}\right)^{k-j}f'\\
    & = & \left(\frac{f'}{f}\right)^{k-2}(f')^2\sum_{j=1}^{k+1}{k+1 \choose j} \\
    & = & \left(\frac{f'}{f}\right)^{k-2}(f')^2(2^{k+1}-1)\\
    & = & \left(\frac{f'}{f}\right)^{k}f^{2}H_{1,1+k}(2).
    \neeq
This completes the proof for the case (a) $m=1$ and any $d>0$.

(b) Assuming (\ref{main}) for $m=1, 2, \ldots k$ and any $d>0$ we shall prove it for $m=k+1$ and any $d>0$.
Since $G_{k+1}(x)=F(x)G_{k}(x)$, by (\ref{der}) and the induction assumption, we obtain
\nbeq
G_{k+1}^{(k+1+d)}(0)
    &  = &
\sum_{j=1}^{k+1+d} {k+1+d \choose j} f^{(j-1)}G_{k}^{(k+1+d-j)}(0)\\
    & = &
  \sum_{j=1}^{d+1} {k+1+d \choose j} f^{(j-1)}G_{k}^{(k+1+d-j)}(0)\\
    & = &
    \sum_{j=1}^{d+1}{k+1+d \choose j}\left(\frac{f'}{f}\right)^{j-2}f'\left(\frac{f'}{f}\right)^{1+d-j}\! \!
    f^{k+1}H_{k,k+1+d-j}(m)\\
    &  = &
    \left(\frac{f'}{f}\right)^{d}f^{k+2} \sum_{j=1}^{k+1+d}{k+1+d \choose j}H_{k,k+1+d-j}(k+1)\\
    &  = &
    \left(\frac{f'}{f}\right)^{d}f^{k+2} \sum_{i=0}^{k+d}{k+1+d \choose i}H_{k,i}(k+1)\\
 &   = &
\left(\frac{f'}{f}\right)^{d}f^{k+2} H_{k+1,k+1+d}(k+2),
    \neeq
where the last equality follows from (\ref{rec_formula}) with $s=k+1$ and $r=k+1+d$. This proves the induction step (b).
Now (iii) follows from (a) and (b).

\vspace{-0.5cm}\section{Concluding Remarks}
\label{sec:4}
We study the distributional
equation $X+T\stackrel{d}{=}Y$, where the shift (translator) $T$ is a sum of i.i.d. random variables without a specified distribution.
The main result in this paper is a characterization of the exponential distribution via a
relationship involving a pair of maxima of i.i.d. continuous random variables.
As a corollary, we prove that the Sukhatme-R\'{e}nyi decomposition of maxima is also a
characterization property for the exponential distribution.

The proof of the main theorem uses a new technique based on an argument from Arnold and Villasenor (2013),
which requires analyticity of the density function. It is an open question if this assumption can be weakened.

\vspace{-0.5cm}\section*{Acknowledgements}
We thank the reviewers and the associate editor for their constructive critique and suggestions. The first author was partially supported by the NFSR at the MES of Bulgaria,
Grant No DFNI-I02/17 while being on leave from the Institute of Mathematics and Informatics at the Bulgarian Academy of Sciences.
\label{sec:5}

\vspace{-0.5cm}\section*{References}

\noindent Arnold, B.C., Balakrishnan N., and Nagaraja, H.N., 2008.  A First Course in

Order Statistics. SIAM, USA, Philadelphia.

\noindent Arnold, B.C. and Villase\~{n}or, J.A., 2013. Exponential characterizations

motivated by the structure of order statistics in samples of size two, Statist

Probab. Lett. 83, 596-601.

\noindent Castano-Martinez, A., Lopez-Blazquez, F., Salamanea-Mino, B., 2012. Random

translations, contractions and dilations of order statistics
and records.

Statistics 46, 57-67.

\noindent Chakraborty, S. and Yanev, G.P., 2013. Characterization of exponential

distribution through equidistribution conditions for consecutive maxima, J.

Statist. Appl. \& Probab. 2, 237-242.

\noindent Conway, J. B., 1978. Functions of One Complex Variable (Graduate Texts in

Mathematics - Vol 11) (v. 1) 2nd Ed., Springer, New York, USA.

\noindent Jovanovic, M., Milosevic, B., Nikitin, Ya. Yu., Obradovic, M., Volkova, K. Yu.,

2015. Tests of exponentiality based on Arnold-Villasenor characterization

and their efficiencies, Comput. Statist. Data Analysis 90, 100-113.

\noindent Ruiz, S.M., 1996. An algebraic identity leading to Wilson's theorem. The Math.

Gazette 80, 579-582.

\noindent Shah, I.A., Khan, A.H., Barakat, H.M., 2014. Random translation, dilation and

contraction of order statistics, Statist. Probab. Lett. 92, 209-214.

\noindent Volkova, K., 2015. Goodness-Of-Fit Tests for Exponentiality Based on Yanev-

Chakraborty Characterization and Their Efficiencies. In: Ed. Nagy, S.,

Proc. 19th European Young Statisticians Meeting, Prague, 156-159.

\noindent Wesolowski, J. and Ahsanullah, M., 2004. Switching order statistics through

random power contractions. Aust. N. Z. J. Statist. 46, 297-303.

\noindent Yanev, G.P. and Chakraborty, S., 2013. Characterizations of exponential

distribution based on sample of size three. Pliska Studia Mathematica

Bulgarica 23, 237-244.

\end{document}